\documentclass[reqno]{amsart}
\usepackage{amssymb,amsmath,amsfonts,amsthm,latexsym}
\usepackage[english]{babel}
\newtheorem{thm}{Theorem}[section]
\newtheorem{cor}[thm]{Corollary}
\newtheorem{lem}[thm]{Lemma}
\newtheorem{prop}[thm]{Proposition}
\newtheorem{defn}[thm]{Definition}
\newtheorem{rem}[thm]{Remark}

\begin{document}
\author{J. R. G\'{o}mez $^{1},$ B.A. Omirov$^2$}

\title{\bf On classification of complex filiform Leibniz algebras}

\maketitle
\begin{abstract}
In this paper we prove that in classifying of complex filiform
Leibniz algebras, for which its naturally graded algebra is
non-Lie algebra, it suffices to consider some special basis
transformations. Moreover, we establish a criterion whether given two
such Leibniz algebras are isomorphic in terms of such
transformations. The classification problem of filiform Leibniz
algebras, for which its naturally graded algebras are non-Lie in
an arbitrary dimension, is reduced to the investigation of the
obtained conditions.
\end{abstract}

\medskip
$^{1}$ Dpto. Matem\'{a}tica Aplicada I, Universidad  de Sevilla,
Avda.  Reina  Mercedes,  s/n. 41012 Sevilla (Spain), e-mail:
\emph{jrgomez@us.es}

$^{2}$ Institute of Mathematics, 29, Do'rmon Yo'li str., 100125, Tashkent (Uzbekistan), e-mail:
\emph{omirovb@mail.ru}

\medskip \textbf{AMS Subject Classifications (2010):
17A32, 17A36, 17B30.}

\textbf{Key words:}  Lie algebra, Leibniz algebra, filiform
Leibniz algebra, natural gradation, classification, adapted basis.

\section{Introduction}

This paper is devoted to the study of Leibniz algebras,
which have been introduced in \cite{Lod}, \cite{Lod-Pir} and
further investigated in many papers, including, for example
\cite{AO1}, \cite{Cas-Lad}-\cite{Frab}. In fact, it is known that
many properties of nilpotent Lie algebras can be extended to the
Leibniz algebras \cite{AO1}, \cite{Barnes},  \cite{Cartan}.

For an arbitrary Leibniz algebra $L$ with a basis $\{e_0, e_1,
\dots, e_n\}$ the table of multiplication is defined by the
products of the basic elements. Namely, the products
$[e_i,e_j]=\sum\limits_{k=0}^{n}\gamma_{ij}^ke_k$ completely
determine products of arbitrary elements of the algebra. The
constants $\gamma_{ij}^k$ are called the structural constants of
the algebra $L$ at the basis $\{e_0, e_1, \dots, e_n\}$.

Thus the problem of classification of algebras can be reduced to
the problem of finding a description of the structural constants
up to a non-degenerate basis transformation. From the Leibniz
identity we have polynomial equalities for the structural
constants:
$$\sum_{l=0}^{n}(\gamma_{jk}^l\gamma_{il}^m -\gamma_{ij}^l\gamma_{lk}^m +\gamma_{ik}^l\gamma_{lj}^m )=0.$$
But the straightforward description of structural constants is
somewhat cumbersome and therefore usually one has to apply
different methods of investigation.

Since the description of all nilpotent Leibniz algebras is
unsolvable task (even in the case of Lie algebras) we reduce our
discussion with restriction on their nilindex. The first step in
this direction was done by M. Vergne in \cite{V}. She classified
naturally graded Lie algebras of maximal nilindex (filiform
algebras) and presented a description of filiform Lie algebras
into sum of naturally graded Lie algebra and its 2-cocycles. We
should note that in the case of Leibniz algebras, unlike the Lie
algebras, the notion of singly-generated algebra have sense (such
nilpotent algebras called zero-filiform Leibniz algebras and
evidently, they have maximal nilindex). In \cite{AO1} the
existence of only one zero-filiform Leibniz algebra in each
dimension was shown and classification of naturally graded
filiform (in Leibniz algebras case they have nilindex equal to the
maximum minus one) is obtained. Also, the description of filiform
Lie algebras were extended to the Leibniz algebras case.

Many authors have studied the classification of nilpotent Lie
algebras for low dimensions. The lists of nilpotent Lie algebras
up to dimension 8 can be found in \cite{G-K} and the
classification of filiform Lie algebras up to dimension 12 can be
obtained from \cite{11-fil_Lie} and \cite{Low}. The extensions of
the classification of filiform Lie algebras of dimension 6 and 7
to the case of Leibniz algebras were obtained in \cite{Akbar} and
\cite{Ikrom}, respectively.

In fact, the classification algorithm of any variety of algebras
with some conditions in fixed dimension consist of the following
four steps:

- finding a basis (an adapted basis) in which the table of
multiplication of an algebra have the most convenient form;

- to reduce the study of all transformations of the adapted basis
to the simple ones;

- to find relations between parameters (structural constants) in
initial and transformed basses;

- present the list of pairwise non-isomorphic algebras such that
any algebra with the considered conditions is isomorphic to an
algebra of the presented list.

A new interesting algorithm for classifying complex filiform Lie
algebras is given in \cite{New method}. However, our algorithm for
the special case studied here is different and enables us to get
newer results.

In the case of filiform Lie algebras, the first two steps of the
algorithm is already obtained \cite{Low}. It the present paper, we
simplify the algorithm of classification for some filiform Leibniz
algebras. In fact, using results of \cite{AO1}, where the families
of filiform Leibniz algebra, for which its naturally graded
algebra is non-Lie, are obtained (i.e. the first step of the
algorithm for such algebras was done), we complete the next two of
the mentioned steps. Therefore, now for the classification of such
filiform Leibniz algebras in an arbitrary finite dimension, we can
start from the analysis of the obtained conditions for structural
constants and present the final list of the algebras. Moreover,
from Theorem \ref{thm44} we can conclude that description of such
algebras in each dimension is an algorithmically solvable problem.

In \cite{Lie-like} some properties of Leibniz filiform algebras, for which its
naturally graded algebra is a Lie algebra were studied.

Throughout the paper the basic field is the field of complex
numbers and in the tables of multiplication we shall omit the
products which are equal to zero.

\section{Preliminaries}

\begin{defn} \label{def21} {\em (\cite{Lod}) A vector space $L$ over a field $F$ with a
multiplication $[-,-] : L \otimes L \to L$ is called} a Leibniz
algebra {\em if it satisfies the following identity:}
$$[x,[y,z]]=[[x,y],z]-[[x,z],y].$$
\end{defn}

Given an arbitrary Leibniz algebra $L,$ we define the lower series
sequence:
$$L^1=L, \ L^{k+1}=[L^{k},L], \ k\geq 1.$$

Now we define the main object of the paper.

\begin{defn} \label{def22} {\em A Leibniz algebra $L$ is said to be} filiform {\em if
$dim L^i=n-i,$ for $2\leq i \leq n$ and $n=dim L.$}
\end{defn}

Note that the notion of filiform Leibniz algebras agrees with the
notion of filiform Lie algebra \cite{V}.

\begin{defn}\label{def23}{\em Given a filiform Leibniz algebra $L,$ put
$L_i=L^i/L^{i+1}, \ 1 \leq i\leq n-1,$ and $gr L = L_1 \oplus
L_2\oplus\dots L_{n-1}.$ Then $[L_i,L_j]\subseteq L_{i+j}$ and we
obtain the graded algebra $gr L$. If $gr L$ and $L$ are
isomorphic, denoted by $gr L=L,$ we say that} the algebra $L$ is
naturally graded.
\end{defn}

In the following theorem, we summarize the results of the works
 \cite{AO1}, \cite{V}.

\begin{thm} \label{thm24} Any complex $(n+1)$-dimensional naturally graded filiform
Leibniz algebra is isomorphic to one of the following pairwise non
isomorphic algebras:

$$\left\{\begin{array}{ll}
[e_0,e_0]=e_{2},&  \\[1mm]
[e_i,e_0]=e_{i+1}, & \  1\leq i \leq {n-1}
\end{array} \right. \quad \left\{\begin{array}{ll}
[e_0,e_0]=e_{2}, &  \\[1mm]
[e_i,e_0]=e_{i+1}, & \  2\leq i \leq {n-1}\\[1mm]
\end{array} \right.$$

$$\left\{\begin{array}{lll} [e_i,e_0]=-[e_0,e_i]=e_{i+1}, &
1\leq i \leq {n-1}\\[1mm]
[e_i,e_{n-i}]=-[e_{n-i},e_i]=\delta (-1)^{i}e_n & 1\leq i\leq n-1.
\end{array} \right.$$
where $\delta \in\{0,1\}$ for odd $n$ and $\delta=0$ for even $n.$
\end{thm}

It should be noted that the first two algebras are non-Lie Leibniz
algebras and the third one is Lie algebra.

Due to the list of Theorem \ref{thm24}, we derive that the set of
all complex filiform Leibniz algebras is decomposed into three
disjoint families of algebras.

\begin{thm}\label{thm25}
An arbitrary complex $(n+1)$-dimensional filiform Leibniz algebra
$L$ is isomorphic to one of the following algebras: \\
$\mu_{1}^{\overline{\alpha}, \theta }=\left\{\begin{array}{ll}
[e_0,e_0]=e_{2},&  \\[1mm]
[e_i,e_0]=e_{i+1}, & \  1\leq i \leq {n-1}\\[1mm]
[e_0,e_1]=\alpha_3e_3 + \alpha_4e_4+...+ \alpha_{n-1}e_{n-1}+ \theta e_n, & \\[1mm]
[e_j,e_1]=\alpha_3e_{j+2} + \alpha_4e_{j+3}+...+
\alpha_{n+1-j}e_n, & \  1\leq j \leq {n-2}
\end{array} \right.$ \\[1mm]

$\mu_{2}^{\overline{\beta}, \gamma }=\left\{\begin{array}{ll}
[e_0,e_0]=e_{2}, &  \\[1mm]
[e_i,e_0]=e_{i+1}, & \  2\leq i \leq {n-1}\\[1mm]
[e_0,e_1]=\beta_3e_3 + \beta_4e_4+...+ \beta_ne_n, & \\[1mm]
[e_1,e_1]=\gamma e_n, & \\[1mm]
[e_j,e_1]=\beta_3e_{j+2} + \beta_4e_{j+3}+...+ \beta_{n+1-j}e_n, &
\ 2\leq j \leq {n-2}
\end{array} \right.$ \\

$ \mu_{3}^{\alpha, \beta, \gamma}=\left\{\begin{array}{lll} [e_0,e_0]=\alpha e_{n}, & \\[1mm]
[e_1,e_1]=\beta e_{n}, & \\[1mm]
[e_i,e_0]=e_{i+1}, &    1\leq i \leq {n-1}\\[1mm]
[e_0,e_1]=-e_2+\gamma e_n, & \\[1mm]
[e_0,e_i]=-e_{i+1}, &    2\leq i \leq {n-1} \\[1mm]
[e_i,e_j]=-[e_j,e_i] \in lin<e_{i+j+1}, e_{i+j+2}, \dots , e_n>, &
1\leq i \leq n-3,\\[1mm]
&  2 \leq j \leq {n-1-i} \\[1mm]
[e_{n-i},e_i]=-[e_i,e_{n-i}]=(-1)^i \delta e_n, & 1 \leq i \leq n-1 \\[1mm]
\end{array} \right.$ \\
where  [ , ] is the multiplication in $L$ and $ \{e_0, e_1, e_2,
.... , e_n \}$ is the basis of the algebra, $\delta \in \{0, 1\}$
for odd $n$ and $\delta=0$ for even $n.$ Moreover, the table of
multiplication of the family $\mu_{3}^{\alpha, \beta, \gamma}$
should satisfy the Leibniz identity.
\end{thm}

\begin{rem} \rm By Theorem \ref{thm25} the first step of the algorithm is done, i.e.
we find the basis in which the table of multiplication of filiform
Leibniz algebra have the most convenient form. It is easy to see
that algebras from $\mu_{1}^{\overline{\alpha}, \theta },$
$\mu_{2}^{\overline{\beta}, \gamma }$ are non-Lie and Lie algebras
belong to the family $\mu_{3}^{\alpha, \beta, \gamma}.$
\end{rem}

\section{On transformations of complex filiform Leibniz algebras.}

Since an arbitrary filiform Leibniz algebra, up to an isomorphism,
belongs to one of the families of Theorem \ref{thm25}, we conclude
that in order to investigate the isomorphisms inside the families,
we need to study the behavior of the parameters (structural
constants) under the action of the non-degenerate change of basis.
Further throughout the paper we shall consider only the first two
families of Theorem \ref{thm25}.

Let $L$ be a complex filiform $(n+1)$-dimensional Leibniz algebra
which is obtained from the naturally graded filiform non-Lie
Leibniz algebras.
\begin{defn}\label{def31}{\em A basis $\{e_0, e_1, \dots, e_n\}$
of an algebra
is said to be} adapted {\em if the multiplication of the algebra
 has the form \ $\mu_1^{\overline{\alpha}, \theta }$ or
$\mu_2^{\overline{\beta}, \gamma }.$}
\end{defn}

Let $L$ be a Leibniz algebra defined on a vector space $V$ and
$\{e_0, e_1, \dots, e_n\}$ is the adapted basis of the algebra $L.$

\begin{defn} \label{def32}{\em A basis transformation $f\in GL(V)$
is said to be} an adapted for the multiplication of the algebra $L$ {\em if
a basis $\{f(e_0), f(e_1), \dots, f(e_n)\}$ is adapted.}
\end{defn}

The closed subgroup of the group $GL(V)$ consisting of adapted
transformations will be denoted by $GL_{ad}(V).$

From the following equalities:
$$\sum\limits_{i=k}^n{a(i)}{\sum\limits_{j=i+p}^n{b(i,j)e_j}}=
\sum\limits_{t=k}^{n-p}{\sum\limits_{i=k}^t{a(i)b(i,t+p)e_{t+p}}=
\sum\limits_{j=k+p}^n{\sum\limits_{i=k}^{j-p}{a(i)b(i,j)e_j}}},$$
we obtain the equality
$$
\sum\limits_{i=k}^n{a(i)}\sum\limits_{j=i+p}^n{b(i,j)e_j}=
\sum\limits_{j=k+p}^n{\sum\limits_{i=k}^{j-p}{a(i)b(i,j)e_j}},
\quad 0\leq p\leq {n-k}, 3\leq k\leq n. \ \eqno(1)$$

\begin{prop}\label{prop33} Let $f\in GL_{ad}(V).$

a) If the algebra $L$ belongs to the family
$\mu_1^{\overline{\alpha}, \theta }$, then $f$ has the following
form:
$$
\left\{ \begin{array}{l}
f(e_0) = \sum_{i=0}^{n}a_ie_i, \\
f(e_1) = (a_0 + a_1)e_1+ \sum_{i=2}^{n-2}a_ie_i +
(a_{n-1}+a_1(\theta - \alpha_n))e_{n-1}+b_ne_n,
\\ {f(e_{i+1})=[f(e_i), f(e_0)], \ 1 \leq i\leq {n-1}} \\
f(e_2)=[f(e_0), f(e_0)].
\end{array} \right.
$$

b) If the algebra $L$ belongs to the family
$\mu_2^{\overline{\beta}, \gamma }$, then $f$ has the following
form:
$$
\left\{ \begin{array}{l}
f(e_0) =\sum_{i=0}^{n}a_ie_i,\\
f(e_1) = b_1e_1 - \frac{a_1b_1\gamma}{a_0}e_{n-1}+b_ne_n
\\ {f(e_{i+1})=[f(e_i), f(e_0)], \ 2 \leq i\leq {n-1}} \\
f(e_2)=[f(e_0), f(e_0)].
\end{array} \right.$$
\end{prop}
\begin{proof} Let $f\in GL_{ad}(V).$ We set
$f(e_0)=\sum\limits_{i=0}^n a_ie_i$ and $f(e_1) =\sum\limits_{j=0}^n b_je_j.$\\

{\bf Case a)}. Consider the product $f(e_2)=[f(e_0), f(e_0)]$.
Using the equality (1) we have
$$
[f(e_0),f(e_0)]=a_0(a_0+a_1)e_2+a_0\sum\limits_{i=3}^n{a_{i-1}e_i}+
a_0a_1(\sum\limits_{i=3}^{n-1}\alpha_{i} e_i+\theta
e_n)+a_{1}^2\sum\limits_{i=3}^n{\alpha_ie_i}+$$
$$a_1\sum\limits_{i=2}^{n-2}{a_i}\sum\limits_
{k=i+2}^n{\alpha_{k+1-i}}e_k=a_0(a_0+a_1)e_2+
a_0\sum\limits_{i=3}^n{a_{i-1}e_i}+
a_1(a_0+a_1)\sum\limits_{i=3}^{n-1}{\alpha_{i} e_i}+a_1(a_0\theta+
a_1\alpha_n)e_n+$$
$$a_1\sum\limits_{i=4}^{n}{a_{i-2}}
\sum\limits_{k=i}^n{\alpha_{k+3-i}}e_k=a_0(a_0+a_1)e_2+
a_0\sum\limits_{i=3}^n{a_{i-1}e_i}+
a_1(a_0+a_1)\sum\limits_{i=3}^{n-1}{\alpha_{i} e_i}+a_1(a_0\theta+
a_1\alpha_n)e_n+$$
$$a_1\sum\limits_{i=4}
^{n}\sum\limits_{i=4}^k{{(a_{i-2}}}{{\alpha_{k+3-i}}e_k})=
a_0(a_0+a_1)e_2+(a_0a_2+a_1(a_0+a_1)\alpha_3)e_3+
\sum\limits_{t=4}^{n-1}(a_0a_{t-1}+a_1(a_0+a_1)\alpha_t+$$
$$a_1\sum\limits_{i=4}^t{a_{i-2}\alpha_{t+3-i}})e_t+(a_0a_{n-1}+
a_1(a_0\theta+a_1\alpha_n)+a_1\sum\limits_{i=4}^n{a_{i-2}
\alpha_{n+3-i}})e_n=f(e_2).$$

Since $[f(e_0),f(e_0)]\in L^2,$ we get $a_0(a_0+a_1)\neq 0.$

Consider the product
$$[f(e_0),f(e_1)]=b_0(a_0+a_1)e_2+\sum\limits_{i=3}^n{c_ie_i}.$$
Since $[f(e_0),f(e_1)]\notin L^2$ and $a_0+a_1\neq 0,$  we
conclude that $b_0=0$.

The properties of the adapted transformation deduce
$f(e_2)=[f(e_1),f(e_0)].$

The product $[f(e_1),f(e_0)]$ in the basis $\{e_0, e_1, \dots,
e_n\}$ has the following form:
$$[f(e_1),f(e_0)]=a_0b_1e_2+(a_0b_2+a_1b_1\alpha_3)e_3+
\sum\limits_{t=4}^{n-1}(a_0b_{t-1}+a_1b_1\alpha_t+
a_1\sum\limits_{i=4}^t{b_{i-2}\alpha_{t+3-i}})e_t+$$
$$(a_0b_{n-1}+a_1b_1\alpha_n+a_1\sum\limits_{i=4}^n
{b_{i-2}\alpha_{n+3-i}})e_n.$$

Comparing the coefficients at the basis elements we get the
conditions to coefficients of the transformation $f$:
\begin{displaymath}
\left\{\begin{array}{l}  a_0+ a_1=b_1, \quad a_2=b_2, \\
a_0a_{t-1}+a_1\sum\limits_{i=4}^t{a_{i-2}\alpha_{t+3-i}}=a_0b_{t-1}+a_1\sum\limits_{i=4}^t{b_{i-2}\alpha_{t+3-i}}, \\
a_0a_{n-1}+a_1(a_0\theta+a_1\alpha_n)+a_1\sum\limits_{i=4}^n{a_{i-2}\alpha_{n+3-i}}=
a_0b_{n-1}+a_1b_1\alpha_n+a_1\sum\limits_{i=4}^n{b_{i-2}\alpha_{n+3-i}}.
\end{array} \right.
\end{displaymath}

From these conditions we have
\begin{displaymath}
\left\{ \begin{array}{l}  b_1=a_0+ a_1, \\
b_i=a_i, \  \ 2\leq i\leq{n-2} \\
b_{n-1}=a_{n-1}+a_1(\theta-\alpha_n).
\end{array} \right.
\end{displaymath}

{\bf Case b)} is proved by a similar way.
\end{proof}

Similarly to \cite{Low},  we introduce  the notion of elementary
transformations for algebras from families
$\mu_1^{\overline{\alpha}, \theta }$ and $\mu_2^{\overline{\beta},
\gamma }.$

\begin{defn}\label{def34} {\em
The following types of the adapted transformations are said to be} elementary: \\
{$\begin{array}{ll} $first type$ \ - \  \tau(a, b, k)= \left\{
\begin{array}{ll}
f(e_0) = e_0 + ae_k  \\
f(e_1) = e_1+ be_k \\
f(e_{i+1})=[f(e_i), f(e_0)], & \ 1 \leq i\leq {n-1}, \ 2 \leq k\leq n  \\
f(e_2)=[f(e_0), f(e_0)]
\end{array} \right. \\ \\

$second type$ \ - \  \vartheta(a, b)= \left\{
\begin{array}{ll}
f(e_0) = ae_0 + be_1  \\
f(e_1) = (a+b)e_1+ b(\theta-\alpha_n)e_{n-1}, & \ a(a+b)\neq 0 \\
f(e_{i+1})=[f(e_i), f(e_0)], & \  1 \leq i\leq {n-1}, \\
f(e_2)=[f(e_0), f(e_0)]
\end{array} \right.\\ \\

$third type$ \ - \ \sigma(b, n)= \left\{
\begin{array}{ll}
f(e_0) = e_0   \\
f(e_1) = e_1+ be_n, \\
f(e_{i+1})=[f(e_i), f(e_0)], & \ 2 \leq i\leq {n-1}, \\
f(e_2)=[f(e_0), f(e_0)]
\end{array} \right. \\ \\

$fourth type$ \ - \  \eta(a, k)= \left\{
\begin{array}{ll}
f(e_0) = e_0 + ae_k  \\
f(e_1) = e_1 \\
f(e_{i+1})=[f(e_i), f(e_0)], & \ 2 \leq i\leq {n-1}, \ 2 \leq k\leq n, \\
f(e_2)=[f(e_0), f(e_0)]
\end{array} \right. 
\end{array}$}\\

{$\begin{array}{ll}

$fifth type$ \ - \  \delta(a, b, d)= \left\{
\begin{array}{ll}
f(e_0) = ae_0 + be_1  \\
f(e_1) = de_1 - \frac{bd \gamma}{a} e_{n-1}, & \ ad\neq 0 \\
f(e_{i+1})=[f(e_i), f(e_0)], & \ 2 \leq i\leq {n-1},  \\
f(e_2)=[f(e_0), f(e_0)]
\end{array} \right.
\end{array}$}\\
where $a, b, d \in \mathbb{C}.$
\end{defn}

Let $f$ be an arbitrary element of the group $GL_{ad}(V),$ then
$f$ can be expressed as superposition of the elementary
transformations.

\begin{prop} \label{prop35}

\

i) Let $f$ has the form a) of Proposition \ref{prop33}. Then
$$f=\tau(a_n, b_n, n)\circ \tau(a_{n-1}, a_{n-1}, n-1)\circ
.... \circ \tau(a_2, a_2, 2)\circ\vartheta(a_0, a_1)$$

ii) Let $f$ have the form b) of Proposition \ref{prop33}. Then
$$f=\sigma(b_n, n)\circ \eta(a_n, n)\circ \eta(a_{n-1}, n-2)\circ ... \circ \eta(a_2, 2)\circ\delta(a_0, a_1, b_1)$$
\end{prop}
\begin{proof} Straightforward.
\end{proof}

For the above decompositions the following is true:

\begin{prop}\label{prop36} \qquad

1) A basis transformation $$g=\tau(a_n, b_n, n)\circ \tau(a_{n-1},
a_{n-1}, n-1)\circ ... \circ \tau(a_2, a_2, 2)$$ does not change
the structural constants of an algebra of the family
$\mu_1^{\overline{\alpha}, \theta }$.

2) A basis transformation $$\varphi=\sigma(b_n, n)\circ \eta(a_n,
n)\circ \eta(a_{n-1}, n-2)\circ ... \circ \eta(a_2, 2)$$ does not
change the structural constants of an algebra of the family
$\mu_2^{\overline{\beta}, \gamma }.$
\end{prop}
\begin{proof} Let us prove the first assertion.

Consider a basis transformation $\tau(a,b,k):$
$$\begin{array}{ll} \tau(a, b, k)= \left\{\begin{array}{ll}
f(e_0) = e_0 + ae_k,  \\
f(e_1) = e_1+ be_k,  \quad \quad \quad  2 \leq k\leq n  \\
f(e_{i+1})=[f(e_i), f(e_0)], \ 1 \leq i\leq {n-1}\\
f(e_2)=[f(e_0), f(e_0)]
\end{array} \right.
\end{array}$$

For $2 \leq k\leq {n-1}$ we put $a=b$ and consider the products
which involve the parameters:
$$[f(e_0),f(e_1)]=
\sum\limits_{i=3}^{n-k+1}{\alpha_i(e_i+ae_{k+i-1})}+
\sum\limits_{i=n-k+2}^{n-1}{\alpha_ie_i}+\theta
e_n=\sum\limits_{i=3}^{n-1}{\alpha_i f(e_i)}+\theta f(e_n),$$ $$
[f(e_1),f(e_1)]=\sum\limits_{i=3}^n{\alpha_ie_i}+a\sum\limits_{i=3}^{n-k+1}{\alpha_ie_{k+i-1}}=
\sum\limits_{i=3}^{n-k+1}{\alpha_i(e_i+ae_{k+i-1})}+
\sum\limits_{i=n-k+2}^n{\alpha_ie_i}=\sum\limits_{i=3}^n{\alpha_i
f(e_i)}.
$$
Therefore, basis transformations $\tau(a, a, k), \ 2\leq
k\leq{n-1}$ for any $a$ do not change the parameters $\alpha_i, \
\theta.$

Analogously, one can check that $\tau(a, b, n)\in GL_{ad}(V)$ does
not change parameters $\alpha_i, \ \theta$ for any value of $a.$

Since a superposition of adapted transformations is again an
adapted transformation, we conclude that transformation
$$g=\tau(a_n, b_n, n)\circ \tau(a_{n-1}, a_{n-1}, n-1)\circ \dots
\circ \tau(a_2, a_2, 2)$$ does not change the structural constants
of family $\mu_1^{\overline{\alpha} \theta }.$

The proof of the second assertion of the proposition is carried
out in a similar way.
\end{proof}

Thus, the problem of the study of all basis transformations is
reduced to the problem of investigation of the second and the
fifth types of elementary transformations for the families
$\mu_1^{\overline{\alpha}, \theta }$ and $\mu_2^{\overline{\beta},
\gamma },$ respectively.

\section{A criterion of isomorphisms of complex filiform non-Lie Leibniz algebras.}

For an arbitrary element $a$ of the Leibniz algebra $L,$ denote
 the operator of right multiplication by $R_a(x)$
 (i.e. $R_a(x)=[x,a]$).

Set
$$R_a^m(x):=[[...[x, \underbrace{ a],a],...,a] }_{m-times} \
\mbox{and} \ R_a^0(x):=x.$$

It should be noted that for an algebra from the first two families
of Theorem \ref{thm25} the following equality holds true:
$$[[e_s,e_1],e_0]=[e_{s+1},e_1], \ 2\leq s\leq n. \eqno(2)$$

Let $L$ be an algebra of the family $\mu_1^{\overline{\alpha},
\theta }$ (respectively, of the family $\mu_2^{\overline{\beta},
\gamma }$ ), then from (2) we derive that for $m\in \mathbb{N}, \
0\leq p\leq n$ (respectively, $0\leq p\leq n, \ p\neq 1$) the
following equality holds:
$$R_{e_1}^m(e_p)=R_{e_0}^{p-1}(R_{e_1}^m(e_0)). \eqno(3)$$

In order to prove the main theorem we need the following lemma.

\begin{lem} \label{lem41} Let $L$ be a  filiform Leibniz algebra
of the first two families from Theorem \ref{thm25}. Then for $2\leq
m\leq \frac{n-1}{2}$ the following equality holds
$$R_{e_1}^m(e_0)=\sum\limits_{i_m=2m+1}^n{\sum\limits_{i_{m-1}=2m+1}^
{i_m}... \sum\limits_{i_1=2m+1}^{i_2}}{\alpha_{i_m+3-i_{m-1}}\cdot
... \cdot \alpha_{i_2+3-i_1}\cdot \alpha_{i_1+3-2(m-1)}e_{i_m}},$$
$\begin{array}{ll} $where$ \ \eta_i =\left\{\begin{array}{ll}
\alpha_i, & $when $L$ belongs to the first family$ \\
\beta_i,  & $when $L$ belongs to the second family.$
\end{array}\right.
\end{array}
$
\end{lem}
\begin{proof}

Let $\eta_i=\alpha_i,$  the case $\eta_i=\beta_i$ is proved
similarly.

We shall use induction by $m$. Using equality (1), for $m=2$ we
have
$$R_{e_1}^2(e_0)=[\sum\limits_{i=3}^{n-1}{\alpha_ie_i+\theta e_n,
e_1}]=\sum\limits_{i=5}^n{\alpha_{i-2}[e_{i-2},
e_1]=\sum\limits_{i=5}^n{\alpha_{i-2}{\sum\limits_{j=i}^n
{\alpha_{j+3-i}e_1}}}}=\sum\limits_{j=5}^n{\sum\limits_{i=5}^j{\alpha_{i-2}\alpha_{j+3-i}e_1}}.$$

Assume that equality of the lemma for $m$ is true. Then the following equalities
$$R_{e_1}^{m+1}(e_0)=[R_{e_1}^m(e_0),e_1]=
[\sum\limits_{i_k=2m+1}^n{\sum\limits_{i_{m-1}=2m+1}^{i_m}...
\sum\limits_{i_1=2m+1}^{i_2}}{\alpha_{i_m+3-i_{k-1}}\cdot ...
\cdot \alpha_{i_2+3-i_1}\cdot \alpha_{i_1+3-(2m+1)}e_{i_m}},
e_1]=$$
$$[\sum\limits_{i_m=2m+3}^n{\sum\limits_{i_{m-1}=2m+3}^{i_m}...
\sum\limits_{i_1=2m+3}^{i_2}}{\alpha_{i_m+3-i_{m-1}}\cdot ...
\cdot \alpha_{i_2+3-i_1}\cdot \alpha_{i_1+3-(2m+3)}e_{i_m-2}},
e_1]=$$
$$\sum\limits_{i_m=2m+3}^n{\sum\limits_{i_{m-1}=2m+3}^{i_m}...
\sum\limits_{i_1=2m+3}^{i_2}}{\alpha_{i_m+3-i_{m-1}}\cdot ...
\cdot\alpha_{i_2+3-i_1}\cdot\alpha_{i_1+3-(2m+3)}\sum\limits_{i_{m+1}=
i_m}^n{\alpha_{i_{m+1}+3-i_m}e_{i_m+1}}}=$$
$$\sum\limits_{i_m=2m+3}^n{{\sum\limits_{i_{m+1}=i_m}^n
\sum\limits_{i_{m-1}=2m+3}^{i_m}...
\sum\limits_{i_1=2m+3}^{i_2}}{\alpha_{i_m+3-i_{m-1}}\cdot ...
\cdot \alpha_{i_2+3-i_1}\cdot \alpha_{i_1+3-(2m+3)}e_{i_m+1}}}=$$
$$\sum\limits_{i_{m+1}=2m+3}^n{{\sum\limits_{i_m=2m+3}^{i_{m+1}}
\sum\limits_{i_{m-1}=2m+3}^{i_m}...
\sum\limits_{i_1=2m+3}^{i_2}}{\alpha_{i_m+3-i_{m-1}}\cdot ...
\cdot \alpha_{i_2+3-i_1}\cdot \alpha_{i_1+3-(2m+3)}e_{i_m+1}}}$$
prove the equality of the lemma for $m+1$ and consequently
completes the proof.
\end{proof}

Since the study of adapted transformations for the family
$\mu_1^{\overline{\alpha}, \theta }$ is reduced to the study of
the following basis transformation
$$\begin{array}{ll} \qquad \qquad \qquad
\left\{\begin{array}{ll}
e_{0}^{'} = \ Ae_0+Be_1, \ \ A(A+B)\neq 0, \\
e_{1}^{'} = \ (A+B)e_1+B(\theta-\alpha_n)e_{n-1},\\
\end{array} \right.
\end{array} \eqno (4)$$
we need the expressions for a new basis. Namely, we have
\begin{cor} \label{cor42}
$$e_{2}^{'}=A(A+B)e_2+B(A+B)\sum\limits_{i=3}^{n-1}{\alpha_i e_i}+
B(A\theta+B\alpha_n),$$

$$e_{k}^{'}=(A+B)(\sum\limits_{i=0}^{k-2}{C_{k-1}^{k-1-i}
A^{k-1-i}B^iR_{e_1}^i(e_{k-i})}+B^{k-1}R_{e_1}^{k-1}(e_0)), \
\eqno(5)$$ where $3\leq k\leq n$ and
$C_{s}^{t}=\frac{s!}{(s-t)!t!}.$
\end{cor}
\begin{proof}

We will prove the corollary by induction on $k$. For $k=2, 3$ we
have
$$e_{2}^{'}=[e_{1}^{'},
e_{0}^{'}]=A(A+B)e_2+B(A+B)\sum\limits_{i=3}^{n-1}{\alpha_i e_i}+
B(A\theta+B\alpha_n),$$

$$e_{3}^{'}=[e_{2}^{'},
e_{0}^{'}]=A^2(A+B)e_3+2AB(A+B)[e_2, e_1] + B^2(A+B)[[e_0,e_1],
e_1]=$$
$$(A+B)(A^2e_3+2ABR_{e_1}(e_2)+ B^2 R_{e_1}^2(e_0)).$$

Suppose that equality (5) is true for $k.$ Taking into account the
equality (2) and the following equalities
$$e_{k+1}^{'}=[e_{k}^{'}, e_0]
=[(A+B)(\sum\limits_{i=0}^{k-2}{C_{k-1}^{k-1-i}A^{p-1-i}
B^iR_{e_1}^{i}(e_{k-i})}+B^{k-1}R_{e_1}^{k-1}(e_0)),Ae_0+Be_1]=$$
$$(A+B)(\sum\limits_{i=0}^{k-2}{C_{k-1}^{k-1-i}A^{k-i}
B^iR_{e_1}^{i}(e_{k+1-i})}+
AB^{k-1}R_{e_1}^{k-1}(e_2)+\sum\limits_{i=0}^{k-2}{C_{k-1}^{k-1-i}
A^{k-1-i}B^{i+1}R_{e_1}^{i+1}(e_{k-i})}+$$
$$B^kR_{e_1}^{k}(e_0))=(A+B)(\sum\limits_{i=0}^{k-2}
{C_{k-1}^{k-1-i}A^{k-i}}B^iR_{e_1}^i(e_{k+1-i})+
AB^{k-1}R_{e_1}^{k-1}(e_2)+\sum\limits_{i=1}^{k-1}{C_{k-1}^{k-i}
A^{k-i}B^{i}R_{e_1}^{i}(e_{k+1-i})}+$$
$$B^kR_{e_1}^{k}(e_0))=
(A+B)(\sum\limits_{i=1}^{k-2}{(C_{k-1}^{k-1-i}+C_{k-1}^{k-i})A^{k-i}
B^iR_{e_1}^{i}(e_{k+1-i}})+
C_{k-1}^{k-1}A^ke_{k+1}+C_{k-1}^{1}AB^{k-1}R_{e_1}^{k-1}(e_2)+$$
$$AB^{k-1}R_{e_1}^{k-1}(e_2)+ B^kR_{e_1}^{k}(e_0))=
(A+B)(\sum\limits_{i=1}^{k-2}{C_{k}^{k-i}A^{k-i}
B^i}R_{e_1}^{i}(e_{k+1-i})+C_{k}^{k}A^ke_{k+1}+
C_{k}^{1}AB^{k-1}R_{e_1}^{k-1}(e_2)+$$
$$B^kR_{e_1}^{k}(e_0))=(A+B)(\sum\limits_{i=0}^{k-1}{C_{k}^{k-i}
A^{k-i}B^iR_{e_1}^{i}(e_{k+1-i})}+B^kR_{e_1}^{k}(e_0))$$ we
complete the proof of the equality (5) for  $k+1.$
\end{proof}

Similarly, for the family  $\mu_2^{\overline{\beta}, \gamma }$
applying the transformation of the type
$$\begin{array}{ll} \left\{\begin{array}{ll}
e_{0}^{'} = \ Ae_0+Be_1, \ \ AD\neq 0, \\
e_{1}^{'} = \ De_1-\frac{BD\gamma}{A}e_{n-1},
\end{array} \right.
\end{array}$$
one can prove the following corollary.
\begin{cor} \label{cor43} For arbitrary $3\leq k\leq
n,$
$$e_{k}^{'}=A(\sum\limits_{i=0}^{k-2}{C_{k-1}^{k-1-i}A^{k-1-i}B^iR_{e_1}^i(e_{k-i})}
+B^{k-1}R_{e_1}^{k-1}(e_0)).$$
\end{cor}
\begin{proof} The proof is carried out in a similar way
as the proof of Corollary \ref{cor42}.\end{proof}

We shall denote an algebra from family $\mu_1^{\overline{\alpha},
\theta }$ (respectively, $\mu_2^{\overline{\beta}, \gamma }$) as
$L(\alpha_3, \alpha_4, \dots, \alpha_n, \theta )$ (respectively,
$L(\beta_3, \beta_4, \dots, \beta_n, \gamma )$).

\begin{thm} \label{thm44}

a) Two algebras  $L( \alpha_3, \alpha_4, \dots, \alpha_n, \theta
)$ and $L^{'}(\alpha'_{3}, \alpha'_{4}, \dots, \alpha'_{n},
\theta')$ are isomorphic if and only if there exist $A, B\in
\mathbb{C}$ such that $A(A+B)\neq 0$ and the following conditions
hold:
$$
\begin{array}{l}
\alpha'_3=\frac{(A+B)}{A^2} \alpha_3 \\
\alpha'_t=\frac{1}{A^{t-1}} \left((A+B)\alpha_t -
\sum\limits_{k=3}^{t-1}( C_{k-1}^{k-2}A^{k-2}B\alpha_{t+2-k}
+\right.
C_{k-1}^{k-3}A^{k-3}B^2\sum\limits_{i_1=k+2}^{t}\alpha_{t+3-i_1}
\alpha_{i_1+1-k} +\\
\  + \ \ \ C_{k-1}^{k-4}A^{k-4}B^3\sum\limits_{i_2=k+3}^{t}
\sum\limits_{i_1=k+3}^{i_2}\alpha_{t+3-i_2}\ \ \cdot \ \
\alpha_{i_2+3-i_1} \ \ \cdot \ \
\alpha_{i_1-k} + \ \ \dots \ \ +  \\
\  + C_{k-1}^{1}AB^{k-2}\sum\limits_{i_{k-3}=2k-2}^{t}
\sum\limits_{i_{k-4}=2k-2}^{i_{k-3}} ...
\sum\limits_{i_1=2k-2}^{i_2} \alpha_{t+3-i_{k-3}}\cdot
\alpha_{i_{k-3}+3-i_{k-4}}\cdot ... \cdot
\alpha_{i_2+3-i_1}\alpha_{i_1+5-2k}+\\
\  + \left. B^{k-1}\sum\limits_{i_{k-2}=2k-1}^{t}
\sum\limits_{i_{k-3}=2k-1}^{i_{k-2}} ... \sum\limits_{i_1=2k-1}^{i_2}
\alpha_{t+3-i_{k-2}}\cdot \alpha_{i_{k-2}+3-i_{k-3}}\cdot ... \cdot
\alpha_{i_2+3-i_1}\cdot \alpha_{i_1+4-2k})\cdot \alpha'_{k}\right),
\end{array}
$$
$$\begin{array}{l}
\theta'=\frac{1}{A^{n-1}}\left(A\theta +B\alpha_n-
\sum\limits_{k=3}^{n-1}(
C_{k-1}^{k-2}A^{k-2}B\alpha_{n+2-k}+\right.
C_{k-1}^{k-3}A^{k-3}B^2\sum\limits_{i_1=k+2}^{n}\alpha_{n+3-i_1}
\alpha_{i_1+1-k} +\\
\  +  \ \ \ C_{k-1}^{k-4}A^{k-4}B^3\sum\limits_{i_2=k+3}^{n}
\sum\limits_{i_1=k+3}^{i_2}\alpha_{n+3-i_2}  \ \  \cdot \ \
\alpha_{i_2+3-i_1} \ \ \cdot \ \
\alpha_{i_1-k}\ \  \ + \ \ \ ... \ \ \ + \\
\  + C_{k-1}^{1}AB^{k-2}\sum\limits_{i_{k-3}=2k-2}^{n}
\sum\limits_{i_{k-4}=2k-2}^{i_{k-3}} \ \ ... \ \
\sum\limits_{i_1=2k-2}^{i_2} \alpha_{n+3-i_{k-3}}
\alpha_{i_{k-3}+3-i_{k-4}} ... \alpha_{i_2+3-i_1}
\alpha_{i_1+5-2k} + \\
\  + \left. B^{k-1}\sum\limits_{i_{k-2}=2k-1}^{n}
\sum\limits_{i_{k-3}=2k-1}^{i_{k-2}} \ \ ... \ \
\sum\limits_{i_1=2k-1}^{i_2} \alpha_{n+3-i_{k-2}}\cdot
\alpha_{i_{k-2}+3-i_{k-3}}\cdot ... \cdot \alpha_{i_2+3-i_1}\cdot
\alpha_{i_1+4-2k})\cdot \alpha'_{k}\right),
\end{array}
$$
where $4\leq t\leq n.$

b) Two algebras $L(\beta_3, \beta_4, \dots, \beta_n, \gamma )$ and
$L^{'}(\beta'_{3}, \beta'_{4}, \dots, \beta'_{n}, \gamma')$ are
isomorphic if and only if there exist $A,B\in \mathbb{C}$ such
that $A(A+B)\neq 0$  and the following conditions hold:
$$\begin{array}{ll}
\gamma'=\frac{D^2}{A^{n}} \gamma, \quad \beta'_3=\frac{D}{A^2} \beta_3, \\
\beta'_t=\frac{1}{A^{t-1}}\left(D\beta_t -
\sum\limits_{k=3}^{t-1}(
C_{k-1}^{k-2}A^{k-2}B\beta_{t+2-k}+C_{k-1}^{k-3}A^{k-3}B^2\sum\limits_{i_1=k+2}^{t}\beta_{t+3-i_1}
\cdot \beta_{i_1+1-k} + \right.\\
 \ \ + \ \ C_{k-1}^{k-4}A^{k-4}B^3\sum\limits_{i_2=k+3}^{t}
\sum\limits_{i_1=k+3}^{i_2}\beta_{t+3-i_2} \ \ \cdot \ \
\beta_{i_2+3-i_1} \ \ \cdot \ \
\beta_{i_1-k} \ \ + \ \ ... \  \ + \\
\ +C_{k-1}^{1}AB^{k-2}\sum\limits_{i_{k-3}=2k-2}^{t}
\sum\limits_{i_{k-4}=2k-2}^{i_{k-3}} ...
\sum\limits_{i_1=2k-2}^{i_2} \beta_{t+3-i_{k-3}}
\beta_{i_{k-3}+3-i_{k-4}} ... \beta_{i_2+3-i_1}\beta_{i_1+5-2k}+
\\
\ \left.+ B^{k-1}\sum\limits_{i_{k-2}=2k-1}^{t}
\sum\limits_{i_{k-3}=2k-1}^{i_{k-2}} ....
\sum\limits_{i_1=2k-1}^{i_2} \beta_{t+3-i_{k-2}}
\beta_{i_{k-2}+3-i_{k-3}} ....
\beta_{i_2+3-i_1}\beta_{i_1+4-2k})\beta'_k\right),
\end{array} $$
where $4\leq t\leq n-1.$
$$\begin{array}{ll}
\beta'_n=\frac{BD\gamma}{A^n}+\frac{1}{A^{n-1}}\left(D\beta_n -
\sum\limits_{k=3}^{n-1}(
C_{k-1}^{k-2}A^{k-2}B\beta_{n+2-k}+C_{k-1}^{k-3}A^{k-3}B^2\sum\limits_{i_1=k+2}^{n}\beta_{n+3-i_1}
\cdot \beta_{i_1+1-k} + \right. \\
\ \ \ + \ \ C_{k-1}^{k-4}A^{k-4}B^3\sum\limits_{i_2=k+3}^{n}
\sum\limits_{i_1=k+3}^{i_2}\beta_{n+3-i_2}\ \ \cdot \ \
\beta_{i_2+3-i_1} \ \ \cdot \ \
\beta_{i_1-k} \ \ \ + \ \ \ ... \ \  \ + \  \\
\ +C_{k-1}^{1}AB^{k-2}\sum\limits_{i_{k-3}=2k-2}^{n}
\sum\limits_{i_{k-4}=2k-2}^{i_{k-3}} ...
\sum\limits_{i_1=2k-2}^{i_2} \beta_{n+3-i_{k-3}}
\beta_{i_{k-3}+3-i_{k-4}} ... \beta_{i_2+3-i_1}\beta_{i_1+5-2k}+
\\
\ \left.+ B^{k-1}\sum\limits_{i_{k-2}=2k-1}^{n}
\sum\limits_{i_{k-3}=2k-1}^{i_{k-2}} ....
\sum\limits_{i_1=2k-1}^{i_2} \beta_{n+3-i_{k-2}}
\beta_{i_{k-2}+3-i_{k-3}} ....
\beta_{i_2+3-i_1}\beta_{i_1+4-2k})\beta'_k\right),
\end{array} $$

\end{thm}
\begin{proof} Consider the class $\mu_1^{\overline{\alpha}, \theta }.$
 Let $\{e_0, e_1,\dots, e_n\}$ be a basis of algebra
 $L(\alpha_3, \alpha_4, \dots, \alpha_n, \theta )$, and
 $\{e'_0, e'_1, \dots, e'_n\}$ be a
basis of the algebra $L^{'}(\alpha_{3}^{'}, \alpha_{4}^{'}, \dots,
\alpha_{n}^{'}, \theta^{'}).$

It is easy to see that in algebra $L(\alpha_3, \alpha_4, \dots,
\alpha_n, \theta )$ the following is true:
$$[[e_0, e_1], e_1]=[[e_1, e_1], e_1].$$

We will consider a change of basis (4).

From Lemma \ref{lem41} and equality (3) we obtain
$$R_{e_1}^m(e_{k-m})=\sum\limits_{i_m=k+m}^n
{\sum\limits_{i_{m-1}=k+m}^{i_m}...
\sum\limits_{i_1=2m+1}^{i_2}}{\alpha_{i_m+3-i_{m-1}}...
\alpha_{i_2+3-i_1} \alpha_{i_1+3-(k+m)}e_{i_m}} \eqno(6)$$
where
$m\leq {n-k}$ è $m\leq k\leq n.$

Now we substitute (6) in the equality (5) and using equalities
(1), (3) we obtain the following:
$$e_{k}^{'}=(A+B)(\sum\limits_{i=0}^{k-2}{C_{k-1}^{k-1-i}A^{k-1-i}
B^iR_{e_1}^i(e_{k-1})}+B^{k-1}R_{e_1}^{k-1}(e_0))=(A+B)(A^{k-1}e_k
+$$
$$C_{k-1}^{k-2}A^{k-2}B\sum\limits_{i=k+1}^{n}{\alpha_{i+2-k}e_i}+
C_{k-1}^{k-3}A^{k-3}B^2\sum\limits_{i=k+2}^{n}\sum\limits_{i_1=k+2}^{i}{\alpha_{i+3-i_1}\cdot
\alpha_{i_1+1-k}e_i} +$$
$$C_{k-1}^{k-4}A^{k-4}B^3\sum\limits_{i=k+3}^{n}\sum\limits_{i_2=k+3}^{i}
\sum\limits_{i_1=k+3}^{i_2}{\alpha_{i+3-i_2}\ \ \cdot \ \
\alpha_{i_2+3-i_1}\ \ \cdot \ \ \alpha_{i_1-k}e_i} \ \ +  \  \ ...
\ \ +\ \ $$
$$C_{k-1}^{1}AB^{k-2}\sum\limits_{i=2k-2}^{n}\sum
\limits_{i_{k-3}=2k-2}^{i}
 \ ... \ \sum\limits_{i_1=2k-2}^{i_2} \alpha_{i+3-i_{k-3}}\cdot \ ... \
\cdot \ \alpha_{i_2+3-i_1}\cdot \alpha_{i_1+5-2k}e_1+$$
$$B^{k-1}\sum\limits_{i=2k-1}^{n}\sum\limits_{i_{k-2}=2k-1}^{i} \ \
 ... \ \ \sum\limits_{i_1=2k-1}^{i_2} \alpha_{i+3-i_{k-2}}\ \ \cdot \
 \ ... \
\cdot \ \ \alpha_{i_2+3-i_1}\alpha_{i_1+4-2k}e_i)\ =\ $$
$$=(A+B)(A^{k-1}e_k + C_{k-1}^{k-2}A^{k-2}B\alpha_3e_{k+1}+
(C_{k-1}^{k-2}A^{k-2}B\alpha_4+$$
$$C_{k-1}^{k-3}A^{k-3}B^2\sum\limits_{i_{1}=k+2}^{k+2}
{\alpha_{k+5-i_1}\cdot \alpha_{i_1+1-k}})e_{k+2}+ \ \ ... \ \  +
(C_{k-1}^{k-2}A^{k-2}B\alpha_{t+2-k}+$$
$$C_{k-1}^{k-3}A^{k-3}B^2\sum\limits_{i_{1}=k+1}^{t}{\alpha_{t+3-i_1}
\alpha_{i_1+1-k}}+
C_{k-1}^{k-4}A^{k-4}B^3\sum\limits_{i_{2}=k+3}^{t}
{\sum\limits_{i_1=k+3}^{i_2}{\alpha_{t+3-i_2}\cdot
\alpha_{i_2+3-i_1}\cdot  \ \alpha_{i_1-k}}} +$$
$$+ ... +
C_{k-1}^{1}AB^{k-2}\sum\limits_{i_{k-3}=2k-2}^{t}
{\sum\limits_{i_{k-4}=2k-2}^{i_{k-3}} ...
{\sum\limits_{i_1=2k-2}^{i_2} {\alpha_{t+3-i_{k-3}}\cdot
\alpha_{i_{k-3}+3-i_{k-4}}\cdot \ ... \ \cdot
\alpha_{i_2+3-i_1}\alpha_{i_1+5-2k}}}}+$$
$$B^{k-1}\sum\limits_{i_{k-2}=2k-1}^{t}
{\sum\limits_{i_{k-3}=2k-1}^{i_{k-2}}
...{\sum\limits_{i_1=2k-1}^{i_2}{\alpha_{t+3-i_{k-2}}\cdot
\alpha_{i_{k-2}+3-i_{k-3}}\cdot \ ... \ \cdot
\alpha_{i_2+3-i_1}\cdot \alpha_{i_1+4-2k}}}})e_t +$$
$$+ \ \ ... \ \  + (C_{k-1}^{k-2}A^{k-2}B\alpha_{n+2-k}+
C_{k-1}^{k-3}A^{k-3}B^2\sum\limits_{i_{1}=k+1}^{n}{\alpha_{n+3-i_1}
\alpha_{i_1+1-k}}\ +\ $$
$$C_{k-1}^{k-4}A^{k-4}B^3\sum\limits_{i_{2}=k+3}^{n}
{\sum\limits_{i_1=k+3}^{i_2}{\alpha_{n+3-i_2} \ \ \cdot \ \
\alpha_{i_2+3-i_1}\ \ \cdot \ \ \alpha_{i_1-k}}} + \ \  ... \ \
+$$
$$C_{k-1}^{1}AB^{k-2}\sum\limits_{i_{k-3}=2k-2}^{n}
{\sum\limits_{i_{k-4}=2k-2}^{i_{k-3}} ...
{\sum\limits_{i_1=2k-2}^{i_2}{\alpha_{n+3-i_{k-3}}\cdot
\alpha_{i_{k-3}+3-i_{k-4}}\cdot ... \cdot
\alpha_{i_2+3-i_1}\cdot \alpha_{i_1+5-2k}}}}+$$
$$B^{k-1}\sum\limits_{i_{k-2}=2k-1}^{n}{\sum\limits_{i_{k-3}=2k-1}^
{i_{k-2}} ...
{\sum\limits_{i_1=2k-1}^{i_2}{\alpha_{n+3-i_{k-2}}\cdot
\alpha_{i_{k-2}+3-i_{k-3}}\cdot ... \cdot
\alpha_{i_2+3-i_1}\alpha_{i_1+4-2k}}}})e_n)=$$
$$=(A+B)\left(A^{k-1}e_k+\sum\limits_{t=k+1}^{n}(C_{k-1}^{k-2}A^{k-2}
B\alpha_{t+2-k}+ C_{k-1}^{k-3}A^{k-3}B^2\sum\limits_{i_1=k+2}^{t}
{\alpha_{t+3-i_1}\cdot \alpha_{i_1+1-k}}+\right.$$
$$C_{k-1}^{k-4}A^{k-4}B^3\sum\limits_{i_2=k+3}^{t}
{\sum\limits_{i_1=k+3}^{i_2}{\alpha_{t+3-i_2} \ \ \ \cdot\ \ \
\alpha_{i_2+3-i_1} \ \ \cdot \ \ \alpha_{i_1-k}}} + \ \ \ ... \ \
\ +$$
$$C_{k-1}^{1}AB^{k-2}\sum\limits_{i_{k-3}=2k-2}^{t}
{\sum\limits_{i_{k-4}=2k-2}^{i_{k-3}} ...
{\sum\limits_{i_1=2k-2}^{i_2} {\alpha_{t+3-i_{k-3}}\cdot
\alpha_{i_{k-3}+3-i_{k-4}}\cdot ... \cdot \alpha_{i_2+3-i_1}\cdot
\alpha_{i_1+5-2k}}}}+$$
$$\left. +B^{k-1}\sum\limits_{i_{k-2}=2k-1}^{t}{\sum\limits_{i_{k-3}=2k-1}^{i_{k-2}}
... {\sum\limits_{i_1=2k-1}^{i_2}{\alpha_{t+3-i_{k-2}}\cdot \
\alpha_{i_{k-2}+3-i_{k-3}}\cdot \ ... \ \cdot
\alpha_{i_2+3-i_1}\cdot \alpha_{i_1+4-2k}}}})e_t\right).$$

Consider the following products in the algebra $L^{'}(\alpha_{3}^{'},
\alpha_{4}^{'}, \dots, \alpha_{n}^{'}, \theta^{'})$:
$$[e_{0}^{'}, e_{1}^{'}]=\sum\limits_{k=3}^{n-1}
{\alpha_{k}^{'}e_{k}^{'}}+
\theta^{'}e_{n}^{'}, \quad \quad [e_{1}^{'}, e_{1}^{'}]=
\sum\limits_{k=3}^n{\alpha_{k}^{'}e_{k}^{'}}.$$

Substituting expression $e'_k,$ obtained above, into the product
$[e'_0,e'_1]$ and using the equality (1) with $p=1$,
we derive the following equalities:
$$[e'_0,e'_1]=\sum\limits_{k=3}^{n-1}
\alpha_{k}^{'}(A+B)\left(A^{k-1}e_k +
\sum\limits_{t=k+1}^{n}(C_{k-1}^{k-2}A^{k-2}B\alpha_{t+2-k}+
C_{k-1}^{k-3}A^{k-3}B^2\sum\limits_{i_1=k+2}^{t}
{\alpha_{t+3-i_1}\cdot\alpha_{i_1+1-k}}+\right.$$
$$C_{k-1}^{k-4}A^{k-4}B^3\sum\limits_{i_2=k+3}^{t}
{\sum\limits_{i_1=k+3}^{i_2}{\alpha_{t+3-i_2}\ \ \ \cdot \ \ \
\alpha_{i_2+3-i_1}\ \ \ \cdot\ \ \ \alpha_{i_1-k}}}\ \ \ + \ \ \
...\ \ \ +\ \ $$
$$+C_{k-1}^{1}AB^{k-2}\sum\limits_{i_{k-3}=2k-2}^{t}
{\sum\limits_{i_{k-4}=2k-2}^{i_{k-3}}...
{\sum\limits_{i_1=2k-2}^{i_2}
{\alpha_{t+3-i_{k-3}}\alpha_{i_{k-3}+3-i_{k-4}}...
\alpha_{i_2+3-i_1}\alpha_{i_1+5-2k}}}}+$$
$$\left. B^{k-1}\sum\limits_{i_{k-2}=2k-1}^{t}
{\sum\limits_{i_{k-3}=2k-1}^ {i_{k-2}} \ \ \ ... \ \ \
{\sum\limits_{i_1=2k-1}^{i_2}{\alpha_{t+3-i_{k-2}}
\alpha_{i_{k-2}+3-i_{k-3}}...\alpha_{i_2+3-i_1}\alpha_{i_1+4-2k}}}})
e_t\right) \ \ +$$
$$\theta^{'}A^{n-1}(A+B)e_n=(A+B)\left(\sum\limits_{k=3}^{n-1}{A^{k-1}
\alpha_{k}^{'}e_k}+\sum\limits_{k=3}^{n-1}{\sum\limits_{t=k+1}^{n}
(C_{k-1}^{k-2}A^{k-2}B\alpha_{t+2-k}}+\right.$$
$${C_{k-1}^{k-3}A^{k-3}B^2\sum\limits_{i_1=k+2}^{t}{\alpha_{t+3-i_1}
\alpha_{i_1+1-k}}+C_{k-1}^{k-4}A^{k-4}B^3\sum\limits_{i_2=k+3}^{t}
{\sum\limits_{i_1=k+3}^{i_2}{\alpha_{t+3-i_2}\alpha_{i_2+3-i_1}
\alpha_{i_1-k}}}}
+$$
$$+...+C_{k-1}^{1}AB^{k-2}\sum\limits_{i_{k-3}=2k-2}^{t}
{\sum\limits_{i_{k-4}=2k-2}^{i_{k-3}}...
{\sum\limits_{i_1=2k-2}^{i_2}
{\alpha_{t+3-i_{k-3}}\alpha_{i_{k-3}+3-i_{k-4}}...
\alpha_{i_2+3-i_1}\alpha_{i_1+5-2k}}}}+$$
$$\left.B^{k-1}\sum\limits_{i_{k-2}=2k-1}^{t}{\sum\limits_{i_{k-3}=2k-1}^
{i_{k-2}}...{\sum\limits_{i_1=2k-1}^{i_2}{\alpha_{t+3-i_{k-2}}
\alpha_{i_{k-2}+3-i_{k-3}}...\alpha_{i_2+3-i_1}\alpha_{i_1+4-2k}}}}\right)
\alpha_k^{'}e_t+$$
$$\theta'A^{n-1}e_n)=(A+B)\left(A^2\alpha'_3e_3+\sum\limits_{t=3}^{n-1}
{A^{t-1}\alpha'_te_t}+\sum\limits_{t=3}^{n-1}{\sum\limits_{k=3}^
{t-1}(C_{k-1}^{k-2}A^{k-2}B\alpha_{t+2-k}}+\right.$$
$$C_{k-1}^{k-3}A^{k-3}B^2\sum\limits_{i_1=k+2}^{t}
{\alpha_{t+3-i_1}\alpha_{i_1+1-k}}+C_{k-1}^{k-4}A^{k-4}B^3
\sum\limits_{i_2=k+3}^{t}{\sum\limits_{i_1=k+3}^{i_2}
{\alpha_{t+3-i_2}\alpha_{i_2+3-i_1}\alpha_{i_1-k}}}
+$$
$$+ ... + C_{k-1}^{1}AB^{k-2}\sum\limits_{i_{k-3}=2k-2}^{t}
{\sum\limits_{i_{k-4}=2k-2}^{i_{k-3}}...{\sum\limits_
{i_1=2k-2}^{i_2}
{\alpha_{t+3-i_{k-3}}\alpha_{i_{k-3}+3-i_{k-4}}...
\alpha_{i_2+3-i_1}\alpha_{i_1+5-2k}}}}+$$
$$\left. B^{k-1}\sum\limits_{i_{k-2}=2k-1}^{t}
{\sum\limits_{i_{k-3}=2k-1}^{i_{k-2}}...
{\sum\limits_{i_1=2k-1}^{i_2}
{\alpha_{t+3-i_{k-2}}\alpha_{i_{k-2}+3-i_{k-3}}...
\alpha_{i_2+3-i_1} \alpha_{i_1+4-2k}}}})\alpha_k^{'}e_t\right)+$$
$$(A+B)\left(\theta'A^{n-1}+\sum\limits_{k=3}^{n-1}(C_{k-1}^
{k-2}A^{k-2}B\alpha_{n+2-k}+C_{k-1}^{k-3}A^{k-3}B^2\sum
\limits_{i_1=k+2}^{n}{\alpha_{n+3-i_1}\alpha_{i_1+1-k}}+\right.$$
$$C_{k-1}^{k-4}A^{k-4}B^3\sum\limits_{i_2=k+3}^{n}
{\sum\limits_{i_1=k+3}^{i_2}{\alpha_{n+3-i_2}\ \ \ \dot\ \ \
\alpha_{i_2+3-i_1} \ \ \ \dot\ \ \  \alpha_{i_1-k}}}\ \ \ + \ \ \
... \ \ \ +$$
$$C_{k-1}^{1}AB^{k-2}\sum\limits_{i_{k-3}=2k-2}^{n}
{\sum\limits_{i_{k-4}=2k-2}^{i_{k-3}}...
{\sum\limits_{i_1=2k-2}^{i_2}{\alpha_{n+3-i_{k-3}}
\alpha_{i_{k-3}+3-i_{k-4}}...\alpha_{i_2+3-i_1}
\alpha_{i_1+5-2k}}}}+$$
$$\left.B^{k-1}\sum\limits_{i_{k-2}=2k-1}^{n}{\sum
\limits_{i_{k-3}=2k-1}^{i_{k-2}}...{\sum\limits_{i_1=2k-1}^
{i_2}{\alpha_{n+3-i_{k-2}}\alpha_{i_{k-2}+3-i_{k-3}}...
\alpha_{i_2+3-i_1}\alpha_{i_1+4-2k}}}})\alpha'_k\right)e_n=$$
$$(A+B)\left(A^2\alpha'_3e_3+\sum\limits_{t=3}^{n-1}
(A^{t-1}\alpha'_te_t+\sum\limits_{k=3}^{t-1}
(C_{k-1}^{k-2}A^{k-2}B\alpha_{t+2-k}+
C_{k-1}^{k-3}A^{k-3}B^2\sum\limits_{i_1=k+2}^{t}{\alpha_{t+3-i_1}
\alpha_{i_1+1-k}}+\right.$$
$$C_{k-1}^{k-4}A^{k-4}B^3\sum\limits_{i_2=k+3}^
{t}{\sum\limits_{i_1=k+3}^{i_2}{\alpha_{t+3-i_2}\ \ \ \dot \ \ \
\alpha_{i_2+3-i_1} \ \ \ \dot \ \ \ \alpha_{i_1-k}}} \ \ \ + \ \ \
... \ \ \ +\ \ \ $$
$$C_{k-1}^{1}AB^{k-2}\sum\limits_{i_{k-3}=2k-2}^{t}
{\sum\limits_{i_{k-4}=2k-2}^{i_{k-3}}...
{\sum\limits_{i_1=2k-2}^{i_2}{\alpha_{t+3-i_{k-3}}
\alpha_{i_{k-3}+3-i_{k-4}}...\alpha_{i_2+3-i_1}
\alpha_{i_1+5-2k}}}}+$$
$$B^{k-1}\sum\limits_{i_{k-2}=2k-1}^{t}{\sum\limits_{i_{k-3}=2k-1}^
{i_{k-2}}...{\sum\limits_{i_1=2k-1}^{i_2}{\alpha_{t+3-i_{k-2}}
\alpha_{i_{k-2}+3-i_{k-3}}...\alpha_{i_2+3-i_1}
\alpha_{i_1+4-2k}}}})\alpha_k^{'}e_t)+$$
$$(\theta^{'}A^{n-1}+\sum\limits_{k=3}^{n-1}
(C_{k-1}^{k-2}A^{k-2}B\alpha_{n+2-k}+C_{k-1}^
{k-3}A^{k-3}B^2\sum\limits_{i_1=k+2}^{n}{\alpha_{n+3-i_1}
\alpha_{i_1+1-k}}+$$
$$C_{k-1}^{k-4}A^{k-4}B^3\sum\limits_{i_2=k+3}^{n}
{\sum\limits_{i_1=k+3}^{i_2}{\alpha_{n+3-i_2}\ \ \ \dot \ \ \
\alpha_{i_2+3-i_1} \ \ \ \dot \ \ \ \alpha_{i_1-k}}} \ \ \ + \ \ \
... \ \ \ +\ \ $$
$$C_{k-1}^{1}AB^{k-2}\sum\limits_{i_{k-3}=2k-2}^{n}
{\sum\limits_{i_{k-4}=2k-2}^{i_{k-3}}...
{\sum\limits_{i_1=2k-2}^{i_2}{\alpha_{n+3-i_{k-3}}
\alpha_{i_{k-3}+3-i_{k-4}}...\alpha_{i_2+3-i_1}
\alpha_{i_1+5-2k}}}}+$$
$$\left.B^{k-1}\sum\limits_{i_{k-2}=2k-1}^{n}
{\sum\limits_{i_{k-3}=2k-1}^{i_{k-2}}...
{\sum\limits_{i_1=2k-1}^{i_2}{\alpha_{n+3-i_{k-2}}
\alpha_{i_{k-2}+3-i_{k-3}}...\alpha_{i_2+3-i_1}
\alpha_{i_1+4-2k}}}) \alpha'_k)}e_n \right).$$

Similar expression for $[e'_1, e'_1]$ can be easily obtained by
substitution in the expression $[e'_0, e'_1]$ instead
$\theta'$ the coefficient $\alpha'_n$, namely
$$[e'_1, e'_1]= (A\ +\ B)\left(A^2\alpha_{3}^{'}e_3\ +\
\sum\limits_{t=4}^{n}(A^{t-1}\alpha_{t}^{'}\ +\
\sum\limits_{k=3}^{t-1}(C_{k-1}^{k-2}A^{k-2}B\alpha_{t+2-k}\right.
\ \ +
$$
$$C_{k-1}^{k-3}A^{k-3}B^2\sum\limits_{i_1=k+2}^{t}{\alpha_{t+3-i_1}
\alpha_{i_1+1-k}}+C_{k-1}^{k-4}A^{k-4}B^3\sum\limits_{i_2=k+3}^{t}
{\sum\limits_{i_1=k+3}^{i_2}{\alpha_{t+3-i_2}\alpha_{i_2+3-i_1}
\alpha_{i_1-k}}} + \ ... \ + $$
$$+C_{k-1}^{1}AB^{k-2}\sum\limits_{i_{k-3}=2k-2}^{t}
{\sum\limits_{i_{k-4}=2k-2}^{i_{k-3}}
... {\sum\limits_{i_1=2k-2}^{i_2} {\alpha_{t+3-i_{k-3}}
\alpha_{i_{k-3}+3-i_{k-4}}...
\alpha_{i_2+3-i_1}\alpha_{i_1+5-2k}}}}+$$
$$\left.+B^{k-1}\sum\limits_{i_{k-2}=2k-1}^{t}
{\sum\limits_{i_{k-3}=2k-1}^{i_{k-2}} ...
{\sum\limits_{i_1=2k-1}^{i_2}{\alpha_{t+3-i_{k-2}}
\alpha_{i_{k-2}+3-i_{k-3}} ...
\alpha_{i_2+3-i_1}\alpha_{i_1+4-2k}}}})\alpha_k^{'})e_t\right).$$

On the other hand, we have
$$[e_{0}^{'}, e_{1}^{'}]=[Ae_0+Be_1,
(A+B)e_1+B(\theta-\alpha_n)e_{n-1}]=
(A+B)^2\sum\limits_{t=3}^{n-1}{\alpha_te_t}+(A+B)
(A\theta+B\alpha_n)e_n,$$
$$[e'_1, e'_1]=[(A+B)e_1+B(\theta-\alpha_n)e_{n-1},
(A+B)e_1+B(\theta-\alpha_n)e_{n-1}]=(A+B)\sum\limits_{t=3}^n
{\alpha_te_t}.$$

Comparing the coefficients of the basis elements $e_t$   and
keeping in mind that the coefficient $A+B$ is different from zero,
we get the restrictions, that were outlined in the first assertion
of the theorem.

Using Corollary \ref{cor43}, the assertion b) of the theorem is
proved by applying similar arguments.\end{proof}

\begin{rem} \rm From Theorem \ref{thm44}, we have that $\alpha'_k$
is a polynomial of the form $P_k(A, B, \alpha_1, \alpha_2, \dots,
\alpha_k, \\ \alpha'_1, \alpha'_2, \dots, \alpha'_{k-1}),$ where
parameters $\alpha_1, \alpha_2, \dots, \alpha_k, \alpha'_1,
\alpha'_2, \dots, \alpha'_{k-1}$ are given and coefficients $A, B$
are unknown, but satisfy the condition $A(A+B)\neq 0.$ And
$\beta'_k$ is also a polynomial of the form $Q_k(A, B , D,
\beta_1, \beta_2, \dots, \beta_k, \beta'_1, \beta'_2, \dots,
\beta'_{k-1})$, where parameters $\beta_1, \beta_2, \dots,
\beta_k, \beta'_1, \beta'_2, \dots, \beta'_{k-1}$ are given and
coefficients $A, B , D$ are unknown, but satisfy the condition
$ AD\neq 0.$ Therefore, the finding of
parameters $\alpha'_k$ and $\beta'_k$ are
recursive procedures. Consequently, we conclude that in any given
dimension the problem of the classification (up to an isomorphism)
of complex filiform Leibniz algebras, which are obtained from the
naturally graded filiform non-Lie algebras, is algorithmically
solvable task.
\end{rem}

\bibliographystyle{amsplain}

\end{document}